\documentclass[12pt]{article}
\usepackage{graphicx}
\usepackage{amssymb,amsmath,amsfonts,amsthm,amscd,mathrsfs}
\usepackage{epstopdf}
\usepackage{amsbsy}
\usepackage{psfrag}

\input xy
\xyoption{all}
\CompileMatrices
\newtheorem{theorem}[equation]{Theorem}
\newtheorem{lemma}[equation]{Lemma}

\newtheorem{corollary}[equation]{Corollary}
\newtheorem{proposition}[equation]{Proposition}
\theoremstyle{definition}
\newtheorem{definition}[equation]{Definition}

\theoremstyle{remark}
\newtheorem{remark}[equation]{Remark}
\newtheorem{remark/definition}[equation]{Remark-Definition}
\newtheorem{example}[equation]{Example}

\numberwithin{equation}{subsection}

\newcommand{\abs}[1]{| #1|}
\newcommand{\inv}[1]{\ensuremath{#1^{-1}}}
\newcommand{\norm}[2]{\ensuremath{\mathrm{N}}_{#1}(#2)}

\newcommand{\A}{\ensuremath{\mathbb{A}}}
\newcommand{\cA}{\ensuremath{\mathcal{A}}}

\newcommand{\N}{\ensuremath{\mathbb{N}}}
\newcommand{\Z}{\ensuremath{\mathbb{Z}}}

\renewcommand{\P}{\ensuremath{\mathbb{P}}}

\newcommand{\ig}{\ensuremath{I[G]}}

\newcommand{\zg}{\ensuremath{\Z[G]}}

\newcommand{\gal}{\ensuremath{\textrm{Gal}}}

\newcommand{\chara}{\ensuremath{\textrm{char}}}

\newcommand{\slashfrac}[2]{\ensuremath{\raise1ex\hbox{#1}\kern-.2em/\kern-.30em\lower1ex\hbox{#2}}}

\newcommand{\Hom}{\ensuremath{\textrm{Hom}}}

\newcommand{\eqnref}[1]{(\ref{#1})}

\newcommand{\atg}{\ensuremath{A_2(G)}}

\newcommand{\barm}{\ensuremath{\overline{m}}}
\newcommand{\barn}{\ensuremath{\overline{n}}}

\newcommand{\mpr}{\ensuremath{M_{\scriptscriptstyle{p[c_2]}}^*(G)}}

\newcommand{\sigm}{\ensuremath{\sigma^{\overline{m}}}}
\newcommand{\sigmp}{\ensuremath{\sigma^{\overline{m}'}}}

\newcommand{\aut}{\ensuremath{\mathrm{Aut}}}

\newcommand{\G}{\ensuremath{\langle \sigma_1\rangle \times \ldots \times \langle \sigma_r\rangle}}
\newcommand{\oF}{\ensuremath{\overline{F}}}
\newcommand{\oD}{\ensuremath{\overline{D}}}

\newcommand{\oE}{\ensuremath{\overline{E}}}
\newcommand{\oL}{\ensuremath{\overline{L}}}
\newcommand{\oDE}{\ensuremath{\overline{D_E}}}

\newcommand{\sign}{\ensuremath{\sigma^{\barn}}}
\newcommand{\signp}{\ensuremath{\sigma^{\overline{n}'}}}
\newcommand{\tors}{\ensuremath{\textrm{Tors }}}
\newcommand{\cF}{\ensuremath{\mathcal{F}}}
\newcommand{\cO}{\ensuremath{\mathcal{O}}}
\newcommand{\ind}{\ensuremath{\mathrm{ind}}}
\newcommand{\sD}{\ensuremath{\mathcal{D}}}
\newcommand{\CH}{\ensuremath{\mathrm{CH}}}
\newcommand{\oG}{\ensuremath{\overline{G}}}

\title{Indecomposable $p$-algebras and Galois subfields in generic abelian crossed products}
\author{Kelly McKinnie}
\begin{document}
\maketitle
\abstract
Let $F$ be a Henselian valued field with $\chara(F)=p$ and $D$ a semi-ramified, ``not strongly degenerate'' $p$-algebra.  We show that all Galois subfields of $D$ are inertial.  Using this as a tool we study generic abelian crossed product $p$-algebras, proving among other things that the noncyclic generic abelian crossed product $p$-algebras defined by non-degenerate matrices are indecomposable $p$-algebras.  To construct examples of these indecomposable $p$-algebras with exponent $p$ and large index we study the relationship between degeneracy in matrices defining abelian crossed products and torsion in $CH^2$ of Severi-Brauer varieties.
\setcounter{section}{-1}
\section{Introduction}

Let $G$ be a noncyclic finite abelian $p$-group.  It is shown in \cite{Saltman-Noncrossed} (3.2) that generic abelian crossed product $p$-algebras defined by the group $G$ and a non-degenerate matrix have the property that all Galois subfields of the algebra have Galois group an image of $G$.  Using a modification of Amitsur's comparison technique this result is used to prove the existence of non-crossed product $p$-algebras in \cite{Saltman-Noncrossed} (3.4) and establishes the fact that noncyclic abelian $p$-groups are rigid (see definitions below).  The main result of section \ref{sec:fiveone} of this paper generalizes \cite{Saltman-Noncrossed} (3.2) to a particular class of valued $p$-algebras.  In particular, we show that for $F$ a Henselian field of characteristic $p$ and $D$ a semi-ramified $p$-algebra with separable residue field which is not strongly degenerate (see definition \ref{deg}) all Galois subfields of $D$ are inertial, and in particular are images of $\gal(\oD/\oF)$.  This is the content of Theorem \ref{inertial}.

Let $A/F$ be a $p$-power index division algebra with center $F$.  For $E/F$ any extension of degree prime to $p$, call $A\otimes_FE$ a \emph{prime to $p$ extension} of $A$.  Prime to $p$ extensions of $p$-power index division algebras have been studied in \cite{Brussel}, \cite{RS}, \cite{mor-seth}, \cite{mck-primetop}, for example, proving that certain properties of a division algebra do or do not hold after a prime to $p$ extension.  The generality of Theorem \ref{inertial}, along with the fact that all of the hypotheses hold after a prime to $p$ extension, allow us to deduce three consequences of Theorem \ref{inertial} which are given in sections \ref{sec:fivethree}, \ref{subsection2} and \ref{sec:fivefive} and we now briefly describe.

A finite group $G$ is said to be \emph{rigid} if there exists a $G$-crossed product $A$ with center $F$ such that $A$ is an $H$-crossed product if and only if $H \cong G$ (\cite{fdda's}).  This result was obtained in \cite{Saltman-Noncrossed} (3.2) using generic $p$-Galois extensions and generic polynomials.  As a first consequence of Theorem \ref{inertial} we show in Corollary \ref{gacpinertial} that the generic abelian crossed products from \cite{Saltman-Noncrossed} (3.2) do not become crossed products with respect to any other group after any prime to $p$ extension.  These abelian crossed products were used in \cite{Saltman-Noncrossed} to prove that noncyclic finite abelian $p$-groups are rigid.  Corollary \ref{gacpinertial} shows that the rigidity property holds for these algebras after any prime to $p$ extension.

The second consequence to Theorem \ref{inertial} is given in Corollary \ref{non-crossed}.  For $k$ an infinite field of characteristic $p$ consider $UD(k,p^n)$, the generic division algebras with exponent equal to index equal to $p^n$, and $UD(k,p^m,p^n)$, the generic division algebras with exponent $p^m$ and index $p^n$.  In Corollary \ref{non-crossed} we show for $k$ an infinite field of characteristic $p$ and $n > m \geq 2$,  $UD(k,p^m,p^n)$ and $UD(k,p^n)$ remain non-crossed products after any prime to $p$ extension.  For $n \geq 3$ the generic division $p$-algebra $UD(k,p^n)$ was originally proven to be a non-crossed product algebra in \cite{Saltman-Noncrossed} (3.4).  Our result about prime to $p$ extensions of a generic division algebra is in the spirit of a result of Rowen and Saltman, \cite{RS} (2.1), in which they prove $UD(k,p^n)$ remains a non-crossed product division algebra after any prime to $p$ extension for any field $k$ with $\chara(k)\ne p$ and $n \geq 3$.  It is also in the spirit of \cite{mor-seth} in which it is shown that, for $n>m\geq 2$ and $\chara(k) \ne p$, $UD(k,p^m,p^n)$ is not a crossed product after any prime to $p$ extension.  

The results in Corollary \ref{gacpinertial} and Corollary \ref{non-crossed} extend previous results on rigidity of abelian $p$-groups and non-crossed product $p$-algebras to include statements regarding prime to $p$ extensions.  Since Theorem \ref{inertial} is proven without the use of generic $p$-Galois extensions, these corollaries reprove the original results in \cite{Saltman-Noncrossed} without the use of generic $p$-Galois extensions.  

The third consequence to Theorem \ref{inertial} is given in Theorem \ref{general-indecomposable}.  We show that the $p$-algebras which satisfy the hypotheses of Theorem \ref{inertial} along with the further condition that they are ``non-degenerate'' are indecomposable $p$-algebras and remain so after any prime to $p$ extension (see also Corollary \ref{cor15}).  This is applied to generic abelian crossed products in Corollary \ref{cor16}, proving that generic abelian crossed product $p$-algebras defined by non-degenerate matrices are indecomposable and remain so after any prime to $p$ extension.  


Finally, in section \ref{ch:ten}, Proposition \ref{prop7}, we study the relationship between degeneracy of matrices defining abelian crossed products and torsion in $CH^2$ of the corresponding Severi-Brauer variety.  The proof of this proposition closely follows \cite{Karpenko2} Proposition 5.3.  The relationship between degeneracy and torsion in $CH^2$ is used in Corollary \ref{cor22} to construct abelian crossed product $p$-algebras which have non-degenerate matrix, exponent $p$ and degree $p^n$ for all $n \geq 2$ and all primes $p \ne 2$.  By Corollary \ref{cor16}, the generic abelian crossed products associated to these abelian crossed product $p$-algebras are indecomposable, exponent $p$ and degree $p^n$ $p$-algebras.  For the prime $p=2$ the situation is slightly more complicated and the end result is the construction of an indecomposable 2-algebra with exponent 2 and index 8.  The first such 2-algebra was constructed in \cite{Rowen}.  These examples are given in \ref{example1} and \ref{example2}.

{\bf Acknowledgments} Much of the work in this paper is part of my doctoral dissertation.  I would like to thank my advisor, David J. Saltman, for countless helpful and inspiring conversations.  I would also like to thank Adrian Wadsworth for his many helpful comments and suggestions upon reading earlier drafts of this work.  In particular I would like to thank him for his help in reducing the hypothesis in Theorem \ref{inertial} from maximally complete to Henselian.

\subsection{Definitions and notations}
In this paper a {\bf $p$-algebra} will refer to a finite dimensional, central simple algebra with center a field of characteristic $p>0$ and $p$-power degree.  Given a field $F$ we denote by $\sD(F)$ the collection of finite dimensional division algebras with center $F$.  A valued division algebra $D$ with valuation $v:D^* \to \Gamma\cup \infty$ and $\Gamma$ a totally ordered abelian group will be denoted $(D,v)\in \sD(F)$.  We denote by $V_D$ the valuation ring of $D$, by $U_D$ the group of units of $D$ and by $M_D$ the unique maximal ideal of $V_D$.  We denote by $\Gamma_{D,v}$ the value group of $D$ and by $\oD_v$ the residue division algebra.  We will always drop the subscript $v$ when there is only one valuation under consideration.  Note that all of these definitions carry through for any sub-division algebra of $D$.  A valued division algebra $(D,v) \in \sD(F)$ is said to be {\bf semi-ramified} if $\oD$ is a field and $[\oD:\oF]=|\Gamma_D:\Gamma_F|=\sqrt{[D:F]}$.  In general $(D,v)$ is said to be {\bf defectless} if $[D:F]=[\oD:\oF]\cdot |\Gamma_D:\Gamma_F|$.

Very often in this paper we will assume a semi-ramified division algebra with separable residue field is (or is not) {\bf strongly degenerate}.  This definition is developed in \cite{mck-primetop} and we recall it briefly here.  Let $(K/F,G,z,u,b)$ be an abelian crossed product with $G=\G$, $u=(u_{ij})_{1\leq i,j\leq r}$, $b=(b_i)_{i=1}^r$, and generating elements $z_i$, $1\leq i \leq r$, satisfying $z_iz_j=u_{ij}z_jz_i$ and $z_i^{n_i}=b_i$ where $n_i=\abs{\sigma_i}$.  For $\barn = (n_1,\ldots,n_r)\in \N^r$ set $\sign=\sigma_1^{n_1}\cdot\ldots\cdot\sigma_r^{n_r}$ and set $z^{\barn}=z_1^{n_1}\cdot\ldots\cdot z_r^{n_r}$.  For $\sigm,\sign \in G$, set $u_{\barm,\barn}=z^{\barm}z^{\barn}(z^{\barm})^{-1}(z^{\barn})^{-1}\in K^*$.
\begin{definition}\begin{enumerate}
\item The matrix $u$ is said to be {\bf degenerate} if there are $\sigm,\sign \in G$ and $a,b \in K^*$ such that $\langle \sigm,\sign\rangle$ is noncyclic and $u_{\barm,\barn} = \sigm(a)a^{-1}\sign(b)b^{-1}$.
\item The matrix $u$ is said to be {\bf strongly degenerate} if there exists an element $\sigm\in G$ with prime order and a set of elements $l,k_1,\ldots, k_r \in K^*$ such that for all $1 \leq i \leq r$, $u_{i,\barm}=\sigm(k_i)k_i^{-1}\sigma_i(l)l^{-1}$, where $i$ stands for the standard basis vector $e_i$.
\end{enumerate}\label{deg}\end{definition}
Let $(D,v) \in \sD(F)$ with $(D,v)$ a semi-ramified division algebra with separable residue field $\oD/\oF$.  By \cite{JW} (1.7), the fundamental group homomorphism $\theta_D:\Gamma_D/\Gamma_F \to \gal(\oD/\oF)$ is an isomorphism.  Recall $\theta_D(v(x)+\Gamma_F)$ is the automorphism of $\oD$ sending $\overline{d}\mapsto \overline{xdx^{-1}}$.  Choose a basis $\sigma=\{\sigma_i\}_{i=1}^r$ of the finite abelian group $\gal(\oD/\oF)$ and elements $\pi_i \in D$ so that $\theta_D(\pi_i+\Gamma_F)=\sigma_i$.  Set $\bar{u}_{ij}$ to be the image in $\oD$ of $\pi_i\pi_j\pi_i^{-1}\pi_j^{-1}$ and define $\bar{u}_{\barm,\barn}$ analogously.
\begin{definition} A semi-ramified, valued division algebra $(D,v)$ with separable residue field $\oD/\oF$ is {\bf degenerate} if the matrix $u$ as constructed above is degenerate.  $(D,v)$ is {\bf strongly degenerate} if the matrix $u$ is strongly degenerate. \label{def22}
\end{definition}
These definitions are shown to be well defined in \cite{mck-primetop} Theorem 2.2 
(also see \cite{Boulagouaz}).  Moreover, in \cite{mck-primetop} Theorem 2.2, strong degeneracy in $(D,v)$ is shown to be equivalent to the existence of a non-trivial $p$-power central homogeneous element in $GD_{\gamma}$ with $\gamma \in \Gamma_D-\Gamma_F$, where $GD=\bigoplus_{\gamma \in \Gamma_D}GD_{\gamma}$ is the associated graded division algebra (see section \ref{sec:fiveone}).  The assumptions of not strongly degenerate and non-degenerate on semi-ramified division algebras will be frequently used in this paper along with the results of \cite{mck-primetop}.

\section{Galois subfields in $p$-algebras}\label{sec:fiveone}
In this section we study Galois subfields and their residue fields in $p$-algebras which are semi-ramified with separable residue field and not strongly degenerate.  The main result is Theorem \ref{inertial}.

\subsection{Pure inseparability in associated graded fields}
Throughout this section we assume the following setup and notation.  For any valued field $(K,w)$, let $GK$ be the associated graded field as in \cite{HW} section 5.  Recall, $GK=\bigoplus_{\gamma \in \Gamma_K}GK_{\gamma}$ where $GK_{\gamma}$ is the quotient $GK^{\gamma}/GK^{>\gamma}$ with $GK^{\gamma}=\{k\in K^*|w(k)\geq \gamma\}\cup \{0\}$ and $GK^{>\gamma}=\{k \in K^*|w(k)>\gamma\}\cup \{0\}$.  Let $GK^h=\cup_{\gamma \in \Gamma_K}GK_{\gamma}$ denote the set of homogeneous elements of $GK$ and let $QGK$ denote the quotient field of the integral domain $GK$.  We adopt the notation that a prime following an element will denote its image in the associated graded ring, that is, given $a \in K$, $a'\in GK$ is the associated homogeneous element in $GK_{w(a)}$.  

Let $(D,v) \in \mathcal{D}(F)$ be a valued division $p$-algebra, semi-ramified with $\oD$ separable over $\oF$.  Let $GD=\bigoplus_{\gamma \in \Gamma_D}GD_{\gamma}$ be the associated graded division algebra defined the same way as $GK$.  In general the center of $GD$ may be strictly larger than $GF$, however, by \cite{Boulagouaz2} Corollary 4.4, under our assumptions on $D$ the center of $GD$ is $GF$.  By \cite{mck-primetop} Theorem 2.2
, if $(D,v)$ is ``not strongly degenerate" then there are no non-trivial $p$-power central homogeneous elements in $GD_{\gamma}$ with $\gamma \in \Gamma_D-\Gamma_F$.  In Theorem \ref{inertial} we will use the not strongly degenerate assumption on Henselian valued $(D,v)$ to conclude that all Galois subfields of $D$ are inertial.  In order to do this we need to explore the connection between ramification in a subfield and purely inseparable totally ramified extensions of $GF$ in $GD$.  This is the topic of this section which culminates with Proposition \ref{p1}, the key ingredient in the proof of Theorem \ref{inertial}.  The results in this section are those of Adrian Wadsworth and were obtained via private communication.
\begin{lemma}
Let $(F,v)$ be a valued field and let $L/F$ be a finite extension of fields.  Assume the valuation $v$ extends uniquely from $F$ to $L$ and $L/F$ is normal.  Then $QGL$ is normal over $QGF$.
\label{l1}
\end{lemma}
\proof To show $QGL$ is normal over $QGF$ we need to show that $QGL$ is the splitting field for a set $\Lambda=\{\lambda(x)\}$ of monic polynomials in $QGF[x]$.  By the proof of \cite{HW} Proposition 2.1, $QGL$ has a $QGF$-basis consisting of homogeneous elements of $GL$.  Let $a'\in GL^h$ be an element of such a basis and choose $a\in L$ with image $a'$ in $GL$.  Let $f_a(x) \in F[x]$ be the minimal polynomial of $a$ over $F$.  Then $f_a(a)=0$ and $f_a(x)$ factors as $f_a(x)=\prod_{j=1}^{m}(x-a_j)$ in $L[x]$ since $L$ is normal (assume $a=a_1$).  For each $1\leq j \leq m$ there exists an $F$-isomorphism $\sigma_j:F(a)\cong F(a_j)$ sending $a$ to $a_j$.  Since $L/F$ is normal these $F$-isomorphisms lift to automorphisms $\sigma_j:L \to L$ satisfying $\sigma_j(a)=a_j$.  Therefore, since the valuation $v$ extends uniquely from $F$ to $L$, we have $v(a)=v(a_j)$ for all $1\leq j \leq m$.  

Let $f_i \in F$ denote the coefficients of $f_a(x)=x^m+f_{1}x^{m-1}+\ldots +f_{m-1}x+f_m$.  By the factorization of $f_a(x)$, $f_i=\pm s_i(a_1,\ldots,a_m)$, where $s_i$ is the $i$-th symmetric polynomial in $m$ variables.  As the $s_i$ are homogeneous polynomials and $v(a_j)=v(a_k)$ for all $1\leq j,k\leq m$, each summand in $f_i$ has the same value.  Therefore, 
\begin{equation}s_i(a_1',\ldots,a_m')=\left\{\begin{array}{ll}s_i(a_1,\ldots,a_m)',&\mathrm{or}\\0&\end{array}\right.\label{e1}\end{equation}
where $a_j'$ is the image of $a_j$ in $GL$ and $s_i(a_1,\ldots,a_m)'$ is the image of $s_i(a_1,\ldots,a_m)$ in $GF$.  Note that $s_m(a_1',\ldots,a_m')=a_1'a_2'\cdot\ldots\cdot a_m'$ is not zero.  Set $\lambda_{a'}(x)=\prod_{j=1}^m(x-a_j')$.  Then $a'=a_1'$ is a root of $\lambda_{a'}(x)$ which factors completely in $QGL[x]$ and by (\ref{e1}) $\lambda_{a'}(x) \in GF[x]\subset QGF[x]$.  Let $\Lambda=\{\lambda_{a'}(x)\}$ with $a'$ running over a $QGF$-basis of $QGL$ consisting of elements in $GL^h$.  $QGL$ is a splitting field for the set of monic polynomials in $\Lambda$ hence $QGL$ is normal over $QGF$.
\endproof

For the rest of this section we will assume $(F,v)$ is a valued field of characteristic $p$ and $L/F$ is a $G$-Galois extension of degree $p^n$.  We will also assume that $v$ extends uniquely from $F$ to $L$, $L/F$ is defectless, and the associated residue field $\overline{L}$ is separable over $\overline{F}$.  Denote the extension of $v$ to $L$ by $v$.  Since $v$ is the only valuation around, all residue fields are with respect to $v$.

Since $v$ extends uniquely from $F$ to $L$ the \emph{decomposition group} of $L/F$, $Z_v(L/F)$, equals $G$.  That is, 
\[Z_v(L/F)=\{\sigma \in G\,|\,v(a)=v(\sigma(a))\, \forall\, a \in F\}=G.\]
Let $T_v(L/F)=\{\sigma \in G\,|\,v(\sigma(a)-a)>0\,\,\forall\, a \in V_L\}$, the \emph{inertia group} of $L/F$.  $T_v(L/F)$ is a normal subgroup of $G$ (\cite{Efrat} Theorem 16.1.1).  Set $E=L^{T_v(L/F)}$, the \emph{inertia subfield} of $L/F$.  By \cite{Efrat} Proposition 16.1.3, $\oL/\oE$, is a purely inseparable extension and therefore $\oL=\oE$ since by assumption $\oL$ is separable over $\oF$.  Furthermore, $\oE/\oF$ is a Galois extension with $\gal(\oE/\oF)\cong \gal(E/F)=G/T_v(L/F)$.  Denote this quotient group by $\oG=\gal(\oE/\oF)$.

Set $[L:E]=p^e$.  Since $L/F$ is defectless, $p^e=|\Gamma_L:\Gamma_F|$.  By \cite{HW} Proposition 2.1, $[QGL:QGF]=[GL:GF]=[GL_0:GF_0]\cdot |\Gamma_L:\Gamma_F|$ and $QGL \cong QGF\otimes_{GF}GL$.  Since $L/F$ is defectless, this implies $[QGL:QGF]=[L:F]=p^n$.  Pictorially we have
\[\xymatrix{
\oL\ar@{=}[d]&\longleftrightarrow&L\ar@{-}[d]^{p^e}&\longleftrightarrow&GL\ar@{-}[d]^{p^e}&\subset&QGL\ar@{-}[d]^{p^e}\\
\oE\ar@{-}[d]^{p^{n-e}}&\longleftrightarrow&E\ar@{-}[d]^{p^{n-e}}&\longleftrightarrow&GE\ar@{-}[d]^{p^{n-e}}&\subset&QGE\ar@{-}[d]^{p^{n-e}}\\
\oF&\longleftrightarrow&F&\longleftrightarrow&GF&\subset&QGF\\
}\]
where the degrees in the rightmost two columns are filled in by the following lemma.
\begin{lemma}
Assume $\chara(F)=p>0$.  The graded field $GL$ is purely inseparable over $GE$ of degree $p^e$.
\label{l2}\end{lemma}
\proof By definition $GL = \bigoplus_{\gamma \in \Gamma_L}GL_{\gamma}$ and $GE=\bigoplus_{\delta \in \Gamma_E}GE_{\delta}$.  Since $E/F$ is the inertia subfield of $L$ and $L/F$ is defectless, $\Gamma_E=\Gamma_F$ and $GL_0=\oL = \oE=GE_0$.  Since each $GL_{\gamma}$ is a one-dimensional vector space over $GL_0=GE_0$ we have 
\[GE=\bigoplus_{\gamma \in \Gamma_F}GE_{\gamma}=\bigoplus_{\gamma \in \Gamma_F}GL_{\gamma}\subset GL.\]
It is now clear that $[GL:GE]=[\Gamma_L:\Gamma_F]=p^e$.  Take $l \in GL_{\gamma}$, a homogeneous element of $GL$.  Then $p^e\gamma \in \Gamma_F$, and therefore $l^{p^e}\in GL_{p^e\gamma}=GE_{p^e\gamma}$.  Since $\chara(F)=p$ this is sufficient to show every element $l \in GL$ satisfies $l^{p^e} \in GE$.
\endproof
Since $QGL$ is the quotient field of $GL$, Lemma \ref{l2} shows $QGL$ is purely inseparable over $QGE$ of degree $[GL:GE]=p^e$.  The field extension $QGE/QGF$ is Galois which can be seen as follows.  Each $\overline{\sigma} \in \oG=\gal(\oE/\oF)$ gives rise to a $GF$-automorphism of $GE$, by acting on the homogeneous pieces of $GE$ which are each a one dimensional vector space over $\oE$.  This association is an injection from $\oG \hookrightarrow \aut_{GF}(GE)$ since it is injective when restricted to the degree 0 piece of $GE$.  Furthermore there is an injection $\aut_{GF}(GE) \hookrightarrow \aut_{QGF}(QGE)$ since an automorphism on a  integral domain extends uniquely to an automorphism on its quotient field.  The composition of these maps is a monomorphism $\oG \hookrightarrow \aut_{QGF}(QGE)$.  Since $|\oG|=[E:F]=[QGE:QGF]$, by Galois Theory $QGE$ is Galois over $QGF$ with Galois group canonically isomorphic to  $\oG$.

Since $QGL$ is normal over $QGF$ by Lemma \ref{l1} and $QGL$ is purely inseparable over $QGE$, we have that the automorphism group 
$\mathrm{Aut}_{QGF}(QGL)$ is isomorphic to $\gal(QGE/QGF)$ and thus also to $\oG$.  Let $QT=QGL^{\oG}$, a purely inseparable extension of $QGF$.  Let $T=QT\cap GL$.  Recall that a graded field $S\supseteq R$ is said to be {\bf totally ramified} over $R$ if $S_0=R_0$.
\begin{proposition}[A. Wadsworth]Let $(F,v)$ be a valued field of characteristic $p>0$ and $L/F$ a $G$-Galois extension of degree $p^n$.  Assume $v$ extends uniquely to $L$, $L/F$ is defectless, and the associated residue field $\overline{L}$ is separable over $\overline{F}$.  Then, $GL$ contains a purely inseparable totally ramified graded field extension of $GF$ of degree $[\Gamma_L:\Gamma_F]$.
\label{p1}
\end{proposition}
\proof Let $T=QT\cap GL$ with $QT$ as above.  Since $v$ extends uniquely to $L$, each $\sigma \in G=\gal(L/F)$ induces a graded automorphism of $GL$ over $GF$ and hence an element of $\aut_{QGF}(QGL)$.  By the comments preceding the proposition regarding the isomorphism $\oG \cong \aut_{QGF}(QGL)$, it is clear that the composite map 
\[\xymatrix{G \ar[r]& \aut_{QGF}(QGL)\ar[r]^(.45){\mathrm{res}}& \gal(QGE/QGF)}\]
is surjective.  Hence $G \to \aut_{QGF}(QGL)$ is surjective.  Thus, each $\tau \in \aut_{QGF}(QGL)$ is induced by the action of a $\sigma \in G$ and $\tau$ restricts to the graded automorphism of $GL$ induced by $\sigma$.  
$T$ is therefore the fixed ring of a family of graded automorphisms of $GL$ and is thus a graded subring of $GL$, and in fact a graded field.  
Moreover, the quotient field of $T$ is $QT$.  This can be seen as follows.  Let $a \in QT$ and write $a=b/c$ for $b,c \in GL$.  Then $a=b'/N(c)$ where $N(c)=\prod_{\overline{\sigma}\in\oG}\overline{\sigma}(c)\in GL$ and $b' \in GL$.  We have $N(c)\in T$ and since $a\in QT$ we also have $b' \in QT\cap GL=T$.  Therefore, $a$ is in the quotient field of $T$ and we have shown the quotient field of $T$ is $QT$.  Thus, $QT\cong QGF\otimes_{GF}T$ and $[T:GF]=[QT:QGF]=[QGL:QGF]/|\oG|=|\Gamma_L:\Gamma_F|$.

Let $p^e=|\Gamma_L:\Gamma_F|$ and take $t \in T$.  Since $t \in QT$, $t^{p^e}\in QGF$, but also $t^{p^e} \in T$, hence $t^{p^e} \in GF$.  Therefore the extension $T$ over $GF$ is purely inseparable.  To see that $T$ is totally ramified over $GF$ note that $T_0$ is a purely inseparable field extension of $GF_0=\overline{F}$ contained in $GL_0=\overline{L}$ and hence we must have $T_0=GF_0$ by the separability assumption on $\oL/\oF$.\endproof

\subsection{Inertial subfields}\label{sec:fivetwo}
Let $(L,w)$ be a valued field and assume $F\subseteq L$ is a subfield with valuation $w|_F$.  $L$ is said to be {\bf inertial} over $F$ if $[L:F]=[\oL:\oF]$ and $\oL$ is separable over $\oF$.  The next theorem, the main one of this section, says that under the appropriate hypotheses all Galois Subfields of a $p$-algebra are inertial.  As shown in Corollary \ref{cor1} this result applies to certain abelian crossed product $p$-algebras defined by not strongly degenerate matrices.  Consequences of Theorem \ref{inertial} are collected in section \ref{consequences}.
\begin{theorem}
Let $(F,v)$ be a Henselian valued field of characteristic $p>0$.  Let $(D,v)\in\mathscr{D}(F)$ be a semi-ramified $p$-algebra with separable residue field $\oD/\oF$ which is not strongly degenerate.  If $L \subset D$ is a subfield, Galois over $F$ with group $G'$, then $L$ is inertial over $F$ and in particular, $G'$ is an image of $G =\textrm{Gal}(\overline{D}/\overline{F})$. \label{inertial}
\end{theorem}
\proof  Let $L/F$ be a Galois subfield of $D$ with group $\gal(L/F)=G'$.  Since $(D,v)$ is a semi-ramified division algebra, it is defectless and any subfield of $D$ is also defectless.  Therefore to show $L$ is inertial it is enough to show $\Gamma_L =\Gamma_F$ since $\oL/\oF$ is separable by the assumption on $\oD/\oF$.  Since $(F,v)$ is Henselian, $v$ extends uniquely to $D$ and thus $L/F$ satisfies the hypotheses of Proposition \ref{p1}.  Therefore $GL$ contains a graded field which is purely inseparable and totally ramified over $GF$  of degree $|\Gamma_L:\Gamma_F|$.  Since $(D,v)$ is not strongly degenerate, by \cite{mck-primetop} Theorem 2.2, 
$GD$ contains no purely inseparable, totally ramified graded field extensions of $GF$.  Therefore $|\Gamma_L:\Gamma_F|=1$ and $L/F$ is an inertial extension.

Since $v$ is Henselian it extends uniquely to $L$ and therefore the decomposition group of $L/F$, $Z_v(L/F)$, is $G'$.  Moreover, since $L/F$ is inertial, the inertia subgroup $T_v(L/F)=G'$.  Thus, $\oL/\oF$ is Galois with $\gal(\oL/\oF)\cong G'$ (\cite{Efrat} 16.1.3).  
Since $\oL\subseteq \oD$, this proves $G'$ is an image of $G=\gal(\oD/\oF)$.  
\endproof
\begin{remark} By \cite{JW} Proposition 1.7, $\gal(\oD/\oF)\cong\Gamma_D/\Gamma_F$, the relative value group.\end{remark}
\begin{remark}Theorem \ref{inertial} is obtained in the author's thesis \cite{M} using the assumption that the center of $D$ is maximally complete.  The proof uses Saltman's generic Galois $p$-extensions and generic polynomials and is much more tedious than the proof given above.  However, it is worth noting that a maximally complete assumption on $(F,v)$ in Theorem \ref{inertial} is enough to obtain all of the results in section \ref{consequences}.
\end{remark}

\section{Corollaries of Theorem \ref{inertial}}\label{consequences}
\subsection{Rigidity in prime to $p$ extensions of crossed products}\label{sec:fivethree}
Let $G=\G$ a finite abelian group.  Let $K/F$ be a $G$-Galois extension of fields, and let $\Delta=(K/F,z_{\sigma},u,b)$ be an abelian crossed product defined by the matrix $u=(u_{ij})\in M_r(K)$ and vector $b=\{b_i\}_{i=1}^r$ (\cite{AS},\cite{mck-primetop}).  Let $F'=F(x_1,\ldots, x_r)$ and $K'=K(x_1,\ldots,x_r)$ with $x_i$ independent indeterminates.  The {\bf generic abelian crossed product} associated to $\Delta$ is defined to be $\cA_{\Delta}=(K'/F',z_{\sigma},u,bx)$ where $bx=\{b_ix_i\}_{i=1}^r$.  In \cite{Saltman-Noncrossed}, Saltman shows that noncyclic abelian $p$-groups are rigid by proving the following theorem.
\begin{theorem}[\cite{Saltman-Noncrossed} Theorem 3.2] Suppose that $K/F$ is a Galois extension of fields of characteristic $p$ with Galois group $G$, a noncyclic finite abelian $p$-group.  Let $\Delta=(K/F,z_{\sigma},u,b)$ be an abelian crossed product with $u$ a non-degenerate matrix.  If $L' \subseteq \mathcal{A}_{\Delta}$ is a subfield, Galois over $F'$, with group $G'$, then $G'$ is an image of $G$.  In particular, if $L'$ is a maximal subfield, $G \cong G'$.
\label{Saltman2}
\end{theorem}
In \cite{AS} 
the existence of abelian crossed products $\Delta$ satisfying the hypotheses of Theorem \ref{Saltman2} is established by proving that abelian crossed products with index equal to exponent are defined by a non-degenerate matrix.  This proves that noncyclic abelian $p$-groups are rigid.  In this section we prove that the examples of Theorem \ref{Saltman2}, which imply that noncyclic abelian $p$-groups are rigid, remain examples after any prime to $p$ extension of the generic abelian crossed product.  That is, the abelian crossed products do not become crossed products with respect to any other group after any prime to $p$ extension.

For any $r\geq 1$ and field $F$ let $F''=F((x_1))\ldots((x_r))$ be the field of iterated Laurent series in $r$ variables.  Let $\Delta=(K/F,z_{\sigma},u,b)$ be an abelian crossed product algebra.  We denote by $\A_{\Delta} = \left(K''/F'', z_{\sigma},u, bx \right)\cong \cA_{\Delta}\otimes_{F'}F''$ the {\bf power series generic abelian crossed product} defined by $\Delta$ (see \cite{mck-primetop}, \cite{Tignol}).  As is well known, $F''$ is Henselian with respect to the standard valuation to $\Z^r$ which is ordered with respect to right-to-left lexicographical ordering (\cite{W} Proposition 3.1).

\begin{corollary}
Let $K/F$ be a $G$-Galois extension of fields with $\chara(F)=p$ and $G=\G$ a noncyclic finite abelian $p$-group.  Let $\Delta=(K/F,z_{\sigma},u,b)$ be an abelian crossed product with $u$ a not strongly degenerate 
matrix.  Let $E''/F''$ be any prime to $p$ extension of $F''$.  Then any Galois subfield of $\mathbb{A}_{\Delta} \otimes_{F''} E''$, has Galois group an image of $G$.\label{cor1}
\end{corollary}
\proof  It is enough to show that $\mathbb{A}_{\Delta}\otimes_{F''}E''$ satisfies the hypotheses of Theorem \ref{inertial}.  Clearly the unique extension of $v$ from $F''$ to $E''$ is Henselian.  Moreover, the $p$-algebra $\mathbb{A}_{\Delta}\otimes_{F''}E''$ is semi-ramified with separable residue field and not strongly degenerate by \cite{mck-primetop} 2.15 and 3.8
.  \endproof

\begin{corollary}
Let $K/F$ be a $G$-Galois extension of fields with $\chara(F)=p$ and $G=\G$ a noncyclic finite abelian $p$-group.  Let $\Delta=(K/F,z_{\sigma},u,b)$ be an abelian crossed product with $u$ a not strongly degenerate 
matrix.  Let $E'/F'$ be any prime to $p$ extension of $F'$.  Then any Galois subfield of $\mathcal{A}_{\Delta} \otimes_{F'} E'$ has Galois group an image of $G$.  \label{gacpinertial}
\end{corollary}
\proof  Let $L'/E'$ be a $G'$-Galois subfield of $\mathcal{A}_{\Delta} \otimes_{F'} E'$.  Let $E''$ be a composite of $E'$ and $F''$ which has degree prime to $p$ over $F''$.  Such a composite exists by, e.g., \cite{mck-primetop} Lemma 3.5
.  Then, $L''=L'\otimes_{E'}E''$ is a Galois subfield of 
\[\cA_{\Delta}\otimes_{F'}E'\otimes_{E'}E'' \cong \cA_{\Delta}\otimes_{F'}F''\otimes_{F''}E'' \cong\A_{\Delta} \otimes_{F''} E''\] 
with group $G'$ (see e.g. \cite{JW} Remark. 5.16(b)).  Therefore, $G'$ is an image of $G$ by Corollary \ref{cor1}. \endproof

\begin{remark}Corollary \ref{gacpinertial} shows that the noncyclic generic abelian crossed product $p$-algebras defined by a not strongly degenerate matrix do not become crossed products with respect to any other group after any prime to $p$ extension.  
\end{remark}
\begin{remark}Since ``the matrix $u$ is not strongly degenerate'' is a weaker condition than ``the matrix $u$ is non-degenerate'', Corollary \ref{gacpinertial} reproves \cite{Saltman-Noncrossed} Theorem 3.2 without the use of generic polynomials and generic Galois $p$-extensions.
\end{remark} 

\subsection{Generic division algebras}\label{subsection2}
Throughout this section let $k$ be an infinite field and $s$ and $r$ positive integers with $s\geq 2$.  Let $UD(k,r)$ be the {\bf generic division algebra} defined over $k$ of index (and exponent) $r$ in $s$ variables (see \cite{LN} chapter 14).  If $\chara(k)\ne p$ and $D=UD(k,p^n)$ with $Z(D)=Z$ then, in \cite{RS} Theorem 2.1, it is proven for $n \geq 3$ and any prime to $p$ extension $E/Z$ that $D \otimes_ZE$ is not a crossed product.   To the contrary, it is shown in \cite{RS} Corollary 1.3 that if $A$ is any division algebra with $\textrm{index}(A)=p^2$ then $A$ becomes a $\Z_p\times\Z_p$-crossed product after a prime to $p$ extension.  Therefore \cite{RS} Corollary 1.3 is the best result possible.    In \cite{mor-seth}, Rowen and Saltman's result is extended to show that if $\textrm{char}(k) \ne p$, then the generic division algebra $UD(k,p^m,p^n)$ with exponent $p^m$ and index $p^n$ has no prime to $p$ extension which is a crossed product provided $n>m\geq 2$.

For the case $\textrm{char}(k)=p$, in \cite{Saltman-Noncrossed} Theorem 3.4, it is shown that $UD(k,p^n)$ is not a crossed product for any $n \geq 3$.  The main result of this section shows, for $\chara(k)=p$ and $n > m \geq 2$,  $UD(k,p^m,p^n)$ and $UD(k,p^n)$ are non-crossed products after any prime to $p$ extension.  The technique used to prove the results in \cite{RS} and \cite{mor-seth}, and that we use here, is a modification of Amitsur's comparison technique.
We state this modification here as we use it in the proof of Corollary \ref{non-crossed}.  Note that it contains no assumption on the characteristic of $k$.
\begin{proposition}[\cite{mor-seth} Proposition 6.1; \cite{RS} Theorem 2.2]  Let $k$ be an infinite field and let $D=UD(k,p^m,p^n)$ with $n\geq m$ and $Z=Z(D)$.  Suppose $D$ has a prime to $p$ extension which is a $G$-crossed product for some group $G$.  If $F$ is any field extension of $k$ and $A/F$ is a central division algebra of index $p^n$ and exponent dividing $p^m$, then $A$ has a prime to $p$ extension which is a $G$ crossed product.
\label{Rowen-Saltman}
\end{proposition}

By Corollary \ref{gacpinertial} we know that there exists abelian crossed product $p$-algebras which do not become crossed products with respect to any other group after any prime to $p$ extension.  In examples \ref{example1} and \ref{example2} from section \ref{ch:ten} we construct such algebras with varying exponent and index.  This allows us to prove the following corollary to Theorem \ref{inertial}.
\begin{corollary}
Let $n>m\geq 2$.  Let $k$ be an infinite field with $\chara(k)=p>0$ and $D=UD(k,p^m,p^n)$ or $D=UD(k,p^n)$ with $Z=Z(D)$.  Then $D \otimes_{Z} E$ is not a crossed product for all prime to $p$ field extensions $E/Z$.\label{non-crossed}
\end{corollary}
\proof  Suppose $E/Z$ is a prime to $p$ extension and $D \otimes_{Z} E$ is a crossed product with group $G$.  By Proposition \ref{Rowen-Saltman} 
it is enough to construct a central division algebra $A/F$ of index $p^n$ and exponent dividing $p^m$ with $F \supseteq k$, so that $A$ is not a $G$-crossed product for all prime to $p$ extensions $E/F$.  For this we use Examples \ref{example1} and \ref{example2} which are developed in section \ref{sec:last}.  Let $G'$ be a noncyclic abelian group of order $p^n$ and $\exp(G')\leq p^m$.  Let $A=\cA_{\Delta(G')}$ be the generic abelian crossed product associated to the abelian crossed product $\Delta(G')$ which is given in \eqnref{eqn76}.  By Corollary \ref{cor16} $\ind(\cA_{\Delta(G')})=|G'|$ and $\exp(\cA_{\Delta(G')})=\exp(G')$.  Moreover, by Corollary \ref{cor22}, the matrix defining $\Delta(G')$ is not strongly degenerate for all primes $p$ since $|G'|\geq p^3$.  Therefore by Corollary \ref{gacpinertial} $A$ does not become a crossed product with respect to any group other than $G'$ after any prime to $p$ extension.  This contradicts our assumption about $D$ since there is more than one non-isomorphic noncyclic abelian group of order $p^n$ and exponent dividing $p^m$ for $n>m\geq 2$.   \endproof

\subsection{Indecomposable $p$-algebras}\label{sec:fivefive}
In this section we use the fact that the Galois subfields of non-degenerate $p$-algebras which satisfy the hypotheses of Theorem \ref{inertial} are inertial to prove that they are also indecomposable and remain indecomposable after any prime to $p$ extension (Theorem \ref{general-indecomposable}).  In particular, we apply this to the generic abelian crossed product $p$-algebras defined by non-degenerate matrices to show that they are indecomposable (Corollary \ref{cor16}).  In order to have meaningful examples of indecomposable $p$-algebras, that is, ones with exponent strictly less than index, we construct in section \ref{ch:ten} 
abelian crossed product $p$-algebras satisfying the hypotheses of Corollary \ref{cor16} which have index $p^n$ and exponent $p$ for $p\ne 2$ and $n \geq 2$. 

Other examples of indecomposable $p$-algebras can be found in \cite{Jacob} where the algebras have exponent $p$ and index $p^2$.  Also see \cite{Karpenko2}, where so called ``generic division algebras'' over fields of arbitrary characteristic and of index $p^n$ and exponent $p$  are proven to be indecomposable for any $n\geq 2$ except if $p=2$ then $n\geq 3$.  Karpenko's examples also remain indecomposable after any prime to $p$ extension.  The results in \cite{Karpenko2} are proven by computing torsion in the second Chow group of Severi-Brauer varieties of generic division algebras.  In section \ref{ch:ten} we use the methods from \cite{Karpenko2} to construct abelian crossed product $p$-algebras defined by non-degenerate matrices with prime exponent.  
\begin{theorem}\label{general-indecomposable}
Let $(F,v)$ be a Henselian valued field of characteristic $p$.  Let $D\in \mathcal{D}(F)$ be a semi-ramified $p$-algebra with separable residue field $\oD/\oF$ which is non-degenerate.  Then $D$ is indecomposable.
\end{theorem}

\proof  Set $\mathrm{index}(D) = p^n$ and assume by way of contradiction that $D$ has a non-trivial decomposition, $D \cong D_1 \otimes_F D_2$, with $\mathrm{index}(D_i) = p^{n_i}$, $n_i \geq 1$ and $n_1+n_2=n$.  By \cite{RS} Proposition 1.1, for $i=1,2$ there exists a prime to $p$ extension $E_i/F$ such that $D_i \otimes_F E_i$ contains a cyclic Galois subfield of degree $p$ over $E_i$.  Let $E=E_1E_2$ be a composite extension of $E_1$ and $E_2$ in an algebraic closure of $F$.  Tensoring $D$ up to $E$ we get, $D\otimes_F E \cong (D_1 \otimes_F E) \otimes_{E}(D_2 \otimes_{F}E)$ and each $D_i\otimes_FE$ contains a cyclic Galois subfield of degree $p$ over $E$.  
Let $L_i \subset D_i \otimes_F E$ denote these two cyclic Galois subfields and set $\gal(L_i/E) = \langle \sigma_i' \rangle$.
\begin{equation}D_E=D\otimes_{F}E \cong \underbrace{\left(D_1 \otimes_{F}E\right)}_{\begin{array}{c}\cup \\L_1\end{array}} \otimes_{E}\underbrace{\left(D_2\otimes_{F}E\right)}_{\begin{array}{c}\cup\\ L_2\end{array}}
\label{equation39}\end{equation}

Let $L \subset D_E$ be isomorphic to the tensor product, $L \cong L_1 \otimes_E L_2$, under the isomorphism given in \eqnref{equation39}.  Since $E/F$ is a prime to $p$ extension, $D_E$ is a division algebra and therefore $L$ is a field with Galois group $\gal(L/E) \cong \langle \sigma_1' \rangle \times \langle \sigma_2' \rangle$.  Choose $\pi_i' \in D_i \otimes_F E$ such that inner automorphism by $\pi_i'$ extends the action of $\sigma_i'$ on $L_i$.  Let $\pi_i \in D_E$ be such that $\pi_1 \mapsto \pi_1'\otimes 1$ and $\pi_2\mapsto 1\otimes \pi_2'$ under the isomorphism given in \eqnref{equation39}.  The elements $\pi_1$ and $\pi_2$ commute in $D_E$ since they commute in an isomorphic image.  Since $(F,v)$ is Henselian there exists a unique extension of $v$ to $E$ which is also Henselian.  Furthermore, by \cite{mck-primetop} Theorem 2.15
, $D_E$ satisfies the hypotheses of Theorem \ref{inertial} and therefore $L$ is an inertial Galois subfield of $D_E$.  Let $\overline{L}\subset \oDE$ denote the residue field of $L$.  The group $\gal(\overline{L}/\overline{E}) \cong \langle \sigma_1' \rangle \times \langle \sigma_2' \rangle$ is a quotient of $G:=\gal(\oDE/\oE)$.  Let $\phi:G \to \gal(\overline{L}/\oE)$ be the surjective homomorphism gotten from restriction and consider the composition
\[\xymatrix@1{
\Gamma_{D_E}/\Gamma_E \ar[r]^(.6){\theta_{D_E}}& G \ar[r]^(.35){\phi}& \gal(\overline{L}/\overline{E}).}\]
Under this map the image of $v(x)+\Gamma_E$ is the automorphism gotten by restriction of inner automorphism by $x$ from $V_{D_E}$ to $V_L$, the valuation rings of $D_E$ and $L$ respectively.  In particular, by our choice of $\pi_i$, $\phi \circ \theta_{D_E}(v(\pi_i)+\Gamma_E) = \sigma_i'$.  Let $\theta_{D_E}(v(\pi_i)+\Gamma_E)=\sigma_i \in G$.  Then $\langle \sigma_1,\sigma_2 \rangle$ is not a cyclic group since there is a homomorphic image in which it is not cyclic.  Let $\{\tau_1,\ldots,\tau_r\}$, $r \geq 2$, be a basis of $G$.  Set $\tau^{\barm_i} = \sigma_i$ for vectors $\barm_i \in \N^r$ where $\tau^{\barm_i} = \tau_1^{m_{i1}}\cdot\ldots\cdot\tau_r^{m_{ir}}$ and choose $\rho_i \in D_E$ such that $\theta_{D_E}(v(\rho_i)+\Gamma_E) = \tau_i$.  Set $v_{ij} = \rho_i\rho_j\inv{\rho_i}\inv{\rho_j}$.  By \cite{mck-primetop} Theorem 2.15 
, $D$ remains non-degenerate after any prime to $p$ extension.  Therefore, the matrix $\overline{v} = (\overline{v}_{ij})$ associated to the algebra $D_E$ is a non-degenerate matrix.  For $i=1,2$ we have $\theta_{D_E}(v(\rho^{\barm_i})+\Gamma_E) = \sigma_i$ and therefore, since $\theta_{D_E}$ is an isomorphism, $v(\pi_i)-v(\rho^{\barm_i}) \in \Gamma_E$.  In particular, there exists an $T_i \in E$ and $a_i \in U_{D_E}$ such that $\pi_i = a_i\rho^{\barm_i}T_i$.  For any $x \in D_E$, let $\psi_{x}:D_E \to D_E$ be inner automorphism by $x$.  Then,
\begin{equation}\label{eqn70}
\begin{split}
1&=\pi_1\pi_2\inv{\pi_1}\inv{\pi_2}\\
&=a_1\rho^{\barm_1}T_1a_2\rho^{\barm_2}T_2\inv{(a_1\rho^{\barm_1}T_1)}\inv{(a_2\rho^{\barm_2}T_2)}\\
&=a_1\psi_{\rho^{\barm_1}}(a_2)v_{\barm_1,\barm_2}\psi_{\rho^{\barm_2}}(a_1^{-1})a_2^{-1}
\end{split}
\end{equation}
Each term in the last line of \eqnref{eqn70} has value zero, therefore we may look at its image in the field $\oDE$ and we see,
\[\overline{1} = \overline{v}_{\barm_1,\barm_2} \frac{\overline{a_1}}{\tau^{\barm_2}(\overline{a_1})}\frac{\tau^{\barm_1}(\overline{a_2})}{\overline{a_2}}.\]
This is a contradiction to our assumption that $\overline{v}$ is a non-degenerate matrix.  Therefore $D$ is indecomposable.  \endproof

\begin{remark} Notice that we made the assumption that the $p$-algebra in Theorem \ref{general-indecomposable} was non-degenerate and not merely ``not strongly degenerate''.  Had we only assumed $D$ was not strongly degenerate we would have known that $D$ had no $p$-power central elements after any prime to $p$ extension by \cite{mck-primetop} Corollary 11
.  This condition would imply the algebra is indecomposable in the case of index $p^n$ for $n=2$ or $n=3$.  For, in these cases, any decomposition would have a factor which is a subalgebra of index $p$.  This subalgebra would become cyclic after a prime to $p$ extension and hence the algebra would have a $p$-power central element, a contradiction.  The condition of non-degeneracy on the algebra has the stronger effect of forcing the subfields in the algebra (and prime to $p$ extensions of the algebra) to behave in a rigid fashion.  This is in the spirit of \cite{AS} Lemma 1.7, where degeneracy forces a decomposition after a restriction of scalars.
\end{remark}

\begin{corollary}Let $(F,v)$ be a Henselian valued field with $\chara(F)=p$.  Let $D\in \mathcal{D}(F)$ be a semi-ramified $p$-algebra with separable residue field $\oD/\oF$ which is non-degenerate.  $D$ is indecomposable after any prime to $p$ extension.\label{cor15}
\end{corollary}
\proof   By \cite{mck-primetop} Theorem 2.15
, for every prime to $p$ extension $E/F$, $D_E$ satisfies the conditions of Theorem \ref{general-indecomposable}.  Therefore $D_E$ is indecomposable.  \endproof

\begin{corollary}
Let $\Delta = (K/F,G,z,u,b)$ be an abelian crossed product with noncyclic finite abelian $p$-group $G$.  Assume $u$ is a non-degenerate matrix and $\chara(F)=p$.  If $\mathcal{A}_{\Delta}$ is the generic abelian crossed product associated to $\Delta$, then $\mathcal{A}_{\Delta}$ is indecomposable with $\ind(\cA_{\Delta})=\abs{G}$ and $\exp(\cA_{\Delta})=\mathrm{LCM}(\exp(G),\exp(\Delta))$.  Moreover, $\mathcal{A}_{\Delta}$ remains indecomposable after any prime to $p$ extension. \label{cor16}
\end{corollary}
\proof  The exponent of $\cA_{\Delta}$ is the same as the exponent of $\A_{\Delta}$ which, by \cite{Tignol} Theorem 2.7, is the least common multiple of $\textrm{exp}(G)$ and $\textrm{exp}(\Delta)$.  Moreover, $\cA_{\Delta}$ is a division algebra (\cite{AS}) and thus clearly has index $\abs{G}$.  To prove that $\mathcal{A}_{\Delta}$ is indecomposable after any prime to $p$ extension, by \cite{mck-primetop} Lemma 3.12
, it is enough to show that $\A_{\Delta}$, the power series generic abelian crossed product, is indecomposable after any prime to $p$ extension.  $\A_{\Delta}$ has center a Henselian valued field of characteristic $p$ and, since $u$ is non-degenerate, $\A_{\Delta}$ is a non-degenerate division algebra (\cite{mck-primetop} Lemma 3.8
).  Therefore we may apply Theorem \ref{general-indecomposable} and conclude that $\A_{\Delta}$ is indecomposable.  By Corollary \ref{cor15} $\A_{\Delta}$ remains indecomposable after any prime to $p$ extension.  \endproof 


\section{The Chow group and non-degeneracy}\label{ch:ten}
In this section we make an observation, Proposition \ref{prop7}, connecting degeneracy in matrices defining abelian crossed products and torsion in the Chow group of the associated Severi-Brauer variety.  Using this observation we construct abelian $G$-crossed product $p$-algebras of exponent $p$ and degree $p^n$, $n \geq 2$, $p \ne 2$ which are defined by non-degenerate matrices in Corollary \ref{cor22}.  By Corollary \ref{cor16} the generic abelian crossed products associated to these algebras are indecomposable of exponent $\mathrm{LCM}(\exp(G),p)=\exp(G)$ and index $p^n$ over a field of characteristic $p$ (see example \ref{example1}).  Taking $G$ to be elementary abelian we get indecomposable $p$-algebras of index $p^n$ and exponent $p$.  For the case $p=2$ we construct abelian crossed product 2-algebras of exponent $2$ and degree $2^n$, $n\geq 3$ which are defined by not strongly degenerate matrices in Corollary \ref{cor22}.  In the case of degree $8$ and exponent $2$, the associated generic abelian crossed product 2-algebra is index 8, exponent 2 and indecomposable since it contains no square central elements (see example \ref{example2}).

\subsection{Torsion in $\CH^2$}
Let $F$ be a field.  Given an central simple $F$-algebra $A$, let $X=SB(A)$ be the Severi-Brauer variety of $A$.  Let $\CH^2(X)$ denote the Chow group of codimension 2 cycles on $X$ modulo rational equivalence.  In this section we connect degeneracy of a matrix $u$ defining an abelian crossed product $\Delta$ of exponent $p$ to torsion in $\CH^2(SB(\Delta))$.  

In this section we adopt the notation of \cite{Karpenko2}.  In particular, given a finite dimensional central simple algebra $A/F$ let $X$ be the Severi-Brauer variety of $A$ and let $\P$ be the projective space $X_{\oF}$ where $\oF$ is an algebraic closure of $F$.  Let $K(X)=K_0(X)$ be the Grothendieck group of $X$.  Let $\xi$ be the class of $\cO_{\P}(-1)$ in $K(\P)$.  As mentioned in \cite{Karpenko2} Theorem 3.1, by \cite{Quillen} Theorem 4.1, the restriction map $K(X) \to K(\P)$ is injective and its image is additively generated by $(\textrm{ind}(A^{\otimes i}))\cdot \xi^i$ ($i \geq 0$).  Let $T^iK(X)$ be the topological filtration of $K(X)$ and let $\Gamma^i K(X)$ be the gamma filtration with $G^*TK(X)$ and $G^*\Gamma K(X)$ representing the associated graded rings to these filtrations (\cite{Karpenko2} Definition 2.6 and 2.7).  The following proposition is an observation about degeneracy in abelian crossed products which follows almost directly from the proof of \cite{Karpenko2} Proposition 5.3.  
\begin{proposition}  Let $\Delta=(K/F,G,z,u,b)$ be an abelian crossed product division algebra with exponent $p$ and noncyclic $G\cong\G$.  Let $X=SB(\Delta)$.  If
\begin{enumerate}
\item $p \ne 2$ and $u$ is degenerate, or
\item $p=2$, $r \geq 3$ and $u$ is strongly degenerate,
\end{enumerate}
then $\CH^2(X)$ is torsion free.\label{prop7}
\end{proposition}
\proof By \cite{Karpenko2} 2.15, and \cite{Karpenko5} 3.1, there exists a surjection
\[\xymatrix{
\tors G^2\Gamma K(X) \ar@{->>}[r] & \tors G^2TK(X) \ar[r]^{\simeq} &\tors \CH^2(X)}.\]
Recall here that $G^2\Gamma K(X) =\Gamma^2K(X)/\Gamma^3K(X)$, $G^2TK(X) = T^2K(X)/T^3K(X)$ and we have the equality $T^2K(X)=\Gamma^2K(X)$.  In general the inclusion $\Gamma^iK(X)\subseteq T^iK(X)$ holds (\cite{Karpenko2} 2.14).  The strategy here will be to show that the generator of $\tors G^2\Gamma K(X)$ maps to zero in $\tors G^2 TK(X)$, that is, it is in $T^3K(X)$.  

Assume $p \ne 2$ and $u$ is degenerate.  By \cite{Karpenko2} Proposition 4.13, since $\Delta$ has prime exponent, the group on the left hand side is cyclic and its generator is represented by the element
\begin{equation}
x=p^n(\xi-1)^2-p^{n-2}(\xi^p-1)^2 \in \Gamma^2K(X)=T^2K(X)\nonumber
\end{equation}
where $\textrm{ind}(\Delta)=p^n$, $n\geq 2$.  This formula is valid only for $p\ne 2$.  Since $u$ is degenerate there exists $\sigm,\sign \in G$ and $a,b \in K^*$ so that $\langle \sigm,\sign \rangle$ is noncyclic and 
\[u_{\barm,\barn}=\frac{\sigm(a)}{a}\frac{\sign(b)}{b}.\]
Recall the matrix $u$ is degenerate in this form if and only if the elements $bz^{\barm}$ and $\inv{a}z^{\barn}$ commute in $\Delta$.  We now do a change of basis to find degeneracy with respect to elements of $G$ of order $p$, we will denote them $\sigmp$ and $\signp$.  Let $\langle \sigm,\sign \rangle = \langle \tau_1\rangle \oplus \langle \tau_2 \rangle$.  Set $\tau_i=(\sigm)^{c_i}(\sign)^{d_i}$ for $i=1,2$.  Set 
\[\left((bz^{\barm})^{c_1}(a^{-1}z^{\barn})^{d_1}\right)^{\frac{|\tau_1|}{p}}=b'z^{\barm'} \hspace{.25in}\textrm{ and } \hspace{.25in} \left((bz^{\barm})^{c_2}(a^{-1}z^{\barn})^{d_2}\right)^{\frac{|\tau_2|}{p}}=a'^{-1}z^{\barn'}\]
Then $\sigmp=\tau_1^{|\tau_1|/p}$ and $\signp=\tau_2^{|\tau_2|/p}$ generate different subgroups of $G$ of order $p$.  Moreover, $b'z^{\barm'}$ and $a'^{-1}z^{\barn'}$ commute and therefore
\[u_{\barm',\barn'}=\frac{\sigmp(a')}{a'}\frac{\signp(b')}{b'}.\]
Set $K'=K^{\langle \sigmp,\signp \rangle}$, then $[K:K']=p^2$.  By the proof of 
\cite{AS} Lemma 1.7, $\Delta'=C_{\Delta}(K')$ is a decomposable division algebra.  Write $\Delta'=\Delta_1\otimes_{K'}\Delta_2$.  Since $\Delta'$ has index $p^2$ we have $\textrm{ind}(\Delta_i)=p$ for $i=1,2$.  Let $X'=SB(\Delta')$.  We can now follow Karpenko's proof of \cite{Karpenko2} Proposition 5.3, precisely.  We include the argument here for completeness.  Consider the element 
\[y=p^2(\xi-1)^2-(\xi^p-1)^2 \in T^2K(X').\]
$y$ is a representative of the generator of $\tors G^2 \Gamma K(X')$ by \cite{Karpenko2} 4.13.  Since $\Delta'$ is a decomposable division algebra of index $p^2$ and exponent $p$, the group $\mathrm{CH}^2(X')$ is torsion free (\cite{Karpenko4} Theorem 1).  Hence $y \in T^3K(X')$.  Taking the transfer of $y$ we get
\[\norm{K'/F}{y}=p^n(\xi-1)^2-p^{n-2}(\xi^p-1)^2 = x \in T^3K(X).\]
Consequently, $\tors \CH^2(X)=0$.

Assume $p=2$, $r \geq 3$ and $u$ is strongly degenerate.  Since $u$ is strongly degenerate there exists $\sigm \in G$ of order $2$ and $l,k_1,\ldots,k_r \in K^*$ so that 
\begin{equation}
u_{i,\barm}=\frac{\sigm(k_i)}{k_i}\frac{l}{\sigma_i(l)} \hspace{.25in}\textrm{ for all }i=1,\ldots,r.\label{sdeg}
\end{equation}
For all $\sigma_i \in G$, set $|\sigma_i|=n_i$, a power of $p=2$.  Choose $i,j$ so that $H=\langle \sigma_i^{n_i/p},\sigma_j^{n_j/p},\sigm\rangle$ is elementary abelian of order 8.  Let $K'=K^H$ and set $\Delta'=C_{\Delta}(K')$.  $\Delta'$ is a decomposable division algebra as can be seen as follows.  Since $u_{i,\barm}$ is of the form \eqnref{sdeg}, by \cite{mck-primetop} Lemma 1.7
, $(lz^{\barm})^2$ is central, so in particular $(lz^{\barm})^2\in K'$.  
Set $\Delta'_1$ to be the subalgebra of $\Delta'$ generated by $K^{\langle \sigma_i^{n_i/p}, \sigma_j^{n_j/p}\rangle}$ and $lz^{\barm}$.  $\Delta'_1$ has center $K^H$ and therefore, by the double centralizer theorem (\cite{LN} Theorem 2.8), $\Delta'\cong \Delta'_1\otimes_{K'}\Delta'_2$ where $\Delta'_2=C_{\Delta'}(\Delta'_1)$.  Since $\Delta'_2$ is degree 4, exponent 2, by \cite{Albert} Theorem 11.2, it is a biquaternion algebra and therefore $\Delta'$ is a triquaternion algebra.

Let $X'=SB(\Delta')$.  Again, we now follow Karpenko's proof of \cite{Karpenko2} Proposition 5.3 precisely and we include the argument here for completeness.  By \cite{Karpenko2} Proposition 4.14, since $\Delta$ has exponent 2 the group $\tors G^2\Gamma K(X)$ is cyclic  with generator represented by the element
\[ x=2^{n-1}(\xi-1)^2-2^{n-3}(\xi^2-1)^2 \in \Gamma^2K(X)=T^2K(X),\]
where $\textrm{ind}(\Delta)=2^n$.  This formula is valid only for $p=2$.  Consider the element 
\[y=2^2(\xi-1)^2-(\xi^2-1)^2 \in T^2K(X').\]
By \cite{Karpenko6} Corollary 3.1, since $\Delta'$ is a decomposable division algebra of index $8$ the group $\CH^2(X')$ is torsion free.  Hence $y \in T^3K(X')$.  Taking the transfer of $y$, we get 
\[\norm{K'/F}{y} = 2^{n-1}(\xi-1)^2-2^{n-3}(\xi^2-1)^2 = x\in T^3K(X).\]
Therefore, $\CH^2(X)$ is torsion free.  \endproof
\begin{remark} The converse of Proposition \ref{prop7}(1) is not true.  The author has been able to construct a decomposable abelian crossed product of degree $p^2$ and exponent $p$, $p\ne 2$, which is defined by a non-degenerate matrix.  By \cite{Karpenko2} Proposition 5.3, this algebra has zero torsion in $\CH^2$ of its Severi-Brauer variety.
\end{remark}

\subsection{Constructing abelian crossed products defined by non-degenerate and not strongly degenerate matrices}\label{sec:last}
Let $G \cong\G$, a noncyclic finite abelian group and let $F$ be a field with a $G$-action
.  We recall here the definition of $\Delta'(G)$, the {\bf generic} crossed product with group $G$ and center defined over $F$, as given in \cite{LN} page 84.  
The $G$-crossed product $\Delta'(G)$ is defined by the maximal subfield $L_0$, the matrix $e(u)$ and the vector $e(b)$ which are given as follows.  Let $I[G]$ be the augmentation ideal of the group ring $\Z[G]$.  $\atg$ is defined to be the $G$-lattice which is the kernel in the short exact sequence 
\begin{equation}
\xymatrix{
0 \ar[r]& \atg \ar[r]^(.3)i & P_2(G)=\bigoplus_{i=1}^r \Z[G] d_i \ar[r]^(.7)j& \ig \ar[r] &0
}\label{eqn89}\end{equation}
where $j(d_i)= \sigma_i-1$.  Set $L_0=F(\atg)=q(F[\atg])$, the field of fractions of the commutative group ring $F[\atg]$ which is a domain.  The $G$-actions of $F$ and $\atg$ extend to a $G$-action on $L_0=F(\atg)$ and, since the $G$-action on $\atg$ is faithful, $L_0/L_0^G$ is a $G$-Galois extension of fields.  Let $e:\atg \to F[\atg]$ be the canonical injection taking the additive group $\atg$ to the multiplicative subgroup of $F[\atg]$ with coefficient 1.  For $1 \leq i,j \leq r$, set 
\begin{eqnarray*}
b_i &=& \norm{i}{d_i}=\left(\sum_{j=0}^{|\sigma_i|-1}\sigma_i^j\right)d_i\in\atg \\
u_{ij} &=& (\sigma_i-1)d_j-(\sigma_j-1)d_i\in\atg.\end{eqnarray*}
Define $e(u)=(e(u_{ij}))\in M_r(L_0^*)$ and $e(b)=\{e(b_i)\}_{i=1}^r\in (L_0^*)^{r}$.  
\begin{definition}
Let $G=\G$ be a noncyclic finite abelian group and $F$ a field with a $G$-action.  Define
\[\Delta'(G)=(L_0/L_0^G,z_{\sigma},e(u),e(b)).\]
$\Delta'(G)$ is the {\bf generic crossed product} associated to the finite abelian group $G$ and the field $F$.  $\Delta'(G)$ exists since the matrix $e(u_{ij})$ and the vector $e(b_i)$ satisfy the necessary conditions outlined in \cite{AS} Theorem 1.3.
\end{definition}
Let $c:G\times G \to L_0$ be the 2-cocycle defining $\Delta'(G)$.  That is $\Delta'(G)=(L_0/L_0^G,G,c)$ and $c(\sigma^{\barm_1},\sigma^{\barm_2})=z^{\barm_1}z^{\barm_2}(z^{\barm_1+\barm_2 \,(\textrm{mod }\barn)})^{-1}$, where $\barm_1+\barm_2 \,(\textrm{mod }\barn)$ is the vector $\barm_1+\barm_2$ with the $i$-th entry taken modulo $n_i=|\sigma_i|$. The class of the 2-cocycle $[c]$ is the canonical one in the following sense.  Let $[c_1]\in H^1(G,\ig)$ be the class of the 1-cocycle which is the image of $1 \in \Z=H^0(G,\Z)$ under the long exact sequence of cohomology applied to $0\to \ig \to \zg \to \Z \to 0$.  Let $[c_2]\in H^2(G,\atg)$ be the image of $[c_1]$ under the long exact sequence of cohomology applied to $0 \to \atg \to P_2(G) \to \ig \to 0$.  Let $[e(c_2)]$ be the image of $[c_2]\in H^2(G,\atg)$ in $H^2(G,L_0)$.
\begin{lemma}Let $G$, $[c]$ and $[c_2]$ be defined as above.  Then,
\begin{enumerate}
\item $[c] = [e(c_2)]$, and 
\item $|[c_2]|=|G|$ in $H^2(G,\atg)$.
\end{enumerate}\label{lemma89}
\end{lemma}
\proof Let $c_1:G \to \ig$ be the map defined by $c_1(g)=g-1$.  Then $c_1$ is a 1-cocycle in the class $[c_1]$ described above.  Define the 1-cochain $\varphi:G \to P_2(G)$ as follows.
\[\varphi(\sigm)=\sum_{i=1}^r\sum_{j=0}^{m_{r-i+1}}\sigma_1^{m_1}\cdots\sigma_{r-i}^{m_{r-i}}\sigma_{r-i+1}^j(d_{r-i+1}).\]
The 1-cochain $\varphi$ satisfies $j(\varphi(\sigm))=\sigm-1$.  Set $c_2=\delta_1(\varphi)$, the 1-coboundary of $c_1$.  Then $c_2$ is a 2-cocycle in the canonical class $[c_2]$ defined above.  Instead of showing that $c(\sigma^{\barm_1},\sigma^{\barm_2})=e(c_2(\sigma^{\barm_1},\sigma^{\barm_2}))$ for all $\sigma^{\barm_1}$, $\sigma^{\barm_2} \in G$ directly we show that the two 2-cocycles $c$ and $e(c_2)$ define the same abelian crossed products.  

Let $(L_0/L_0^G,G,e(c_2))=\bigoplus_{g \in G}L_0\cdot w_g$ with $w_gw_h=e(c_2)(g,h)w_{gh}$.  This crossed product is an abelian crossed product, thus $(L_0/L_0^G,G,e(c_2))=(L_0/L_0^G,w_{\sigma},v,d)$ with the second algebra defined by $w_i=w_{\sigma_i}$, $v_{ij}=w_iw_jw_i^{-1}w_j^{-1}$ and $d_i=(w_i)^{n_i}$, where $n_i=|\sigma_i|$.  To prove part 1 of the lemma it is enough to show that $v_{ij}=e(u_{ij})$ and $d_i=e(b_i)$.  It is easy to check that for $i<j$ we have $e(c_2(\sigma_i,\sigma_j))=1$ and for $p+q<n_i$ we have $e(c_2(\sigma_i^p,\sigma_i^q))=1$.  Assume $i<j$.
\begin{eqnarray*}
v_{ij}&=&w_iw_jw_i^{-1}w_j^{-1} = w_iw_j(w_{\sigma_i\sigma_j})^{-1}(w_jw_i(w_{\sigma_j\sigma_i})^{-1})^{-1}\\
&=&e(c_2(\sigma_i,\sigma_j))(e(c_2(\sigma_j,\sigma_i)))^{-1}\\
&=&1\cdot e(\sigma_jd_i-(\sigma_id_j+d_i)+d_j)^{-1}\\
&=&e((\sigma_i-1)d_j-(\sigma_j-1)d_i)\\
&=&e(u_{ij})\\
d_i &=& (w_i)^{n_i}\\
&=&e(c_2(\sigma_i,\sigma_i^{n_i-1}))=e(b_i).
\end{eqnarray*}
This shows $(v_{ij})=(e(u_{ij}))$ for all $i,j$ since $v_{ij}=-v_{ji}$ and $e(u_{ij})=-e(u_{ji})$.  We now prove the second part of the lemma.  Since $P_2(G)$ is a free $G$-module $H^1(G,P_2(G))=0$ (\cite{LN} 12.3).  Therefore we have the injection $0 \to H^1(G,\ig) \to H^2(G,\atg)$ and to show $|[c_2]|=|G|$ it is enough to show $|G|=|[c_1]|$.  From the long exact sequence
\[\ldots \to \zg^G \to \Z^G=\Z \to H^1(G,\ig) \to 0=H^1(G,\zg)\to\ldots\]
we see that $|H^1(G,\ig)|=|\textrm{Im}(\zg^G \to \Z)|=|G|$.  Since $H^1(G,\ig)$ is generated by $[c_1]$ we are done.
\endproof

As seen in the next lemmas the crossed product $\Delta'(G)$ has the nice property of having exponent equal to index equal to $|G|$ when the field $F$ has a nice enough $G$-action.  Let $M$ be a $G$-module.  We say $M$ is an {\bf $H^1$-trivial} $G$-module if $H^1(H,M)=0$ for all subgroups $H\leq G$.

\begin{lemma}  $\atg$ is $H^1$-trivial.
\end{lemma}
\proof  Apply the long exact sequence of Tate cohomology (\cite{Brown} p. 134) to the short exact sequence \eqnref{eqn89}.
\begin{equation}\cdots \to \widehat{H}^0(H,\ig) \to \widehat{H}^1(H,\atg) \to \widehat{H}^1(H,P_2(G)) \to \cdots\label{eqn889}\end{equation}
From \cite{Brown} page 134, $\widehat{H}^0(H,\ig)=\ig^H/(N_H\cdot \ig)$ where $N_H=\sum_{h\in H}h$ is the $H$-norm.  Let $\{g_i\}$ be a set of coset representatives of $G/H$.  Write each $x \in \ig^H\subset \zg$ as $\sum_{i=1}^{|G/H|}\sum_{h \in H}a_{i,h}g_ih$.  Since $xh=x$ for all $h \in H$, $a_{i,h}=a_{i,h'}$ for all $h,h' \in H$.  Therefore, $x\in N_H\cdot \ig$ and in particular, $\widehat{H}^0(H,\ig)=0$ for all $H\leq G$.  By definition $\widehat{H}^i(H,P_2(G))=H^i(H,P_2(G))$ for all $i\geq1$ and $H\leq G$.  Since $P_2(G)$ is a free $\zg$-module it is $H^1$-trivial (\cite{LN} Lemma 12.3).  Combining these two facts and \eqnref{eqn889} we see that $\widehat{H}^1(H,\atg)=H^1(H,\atg)=0$ for all $H\leq G$.  \endproof

\begin{lemma}  Let $G=\G$ a noncyclic finite abelian group and let $F$ be a field with a $G$-action such that $F^*$ is $H^1$-trivial as a $G$-module.  Let $\Delta'(G)$ be defined by $G$ and $F$.  Then, $\exp(\Delta'(G))=\mathrm{ind}(\Delta'(G))=|G|$.  \label{indexp}
\end{lemma}
\proof Let $\Delta'(G)=(L_0/L_0^G,G,c)$.  By Lemma \ref{lemma89}(1), $[c]=[e(c_2)]$, therefore it is enough to show that $e(c_2)$ has order $|G|$ in $H^2(G,L_0)$.  Since $F^*$ and $\atg$ are $H^1$-trivial the direct summand $F^* \oplus \atg$ is an $H^1$-trivial $G$-module.  Therefore, by \cite{LN} Theorem 12.4(a), the map $H^2(G,F^*\oplus\atg) \to H^2(G,L_0)$ which is induced by $e:\atg\to L_0$ is an injection and it is enough to show that $[c_2]=|G|$ in $H^2(G,\atg)$.  This is done in Lemma \ref{lemma89}(2).  Therefore $\textrm{exp}(\Delta'(G))=|G|=\textrm{ind}(\Delta'(G))$.
\endproof

\begin{remark} Let $G=\G$ be a noncyclic $p$-group and let $F$ be a field with $\chara(F)=p$.  Consider $F^*$ as a $G$-module with trivial action by $G$.  Then $F^*$ is $H^1$-trivial since $H^1(H,F^*)=\Hom(H,F^*)=0$.  The last equality follows because $F^*$ has no nontrivial $p$-th roots of unity.  By Lemma \ref{indexp} the corresponding algebra $\Delta'(G)$ has exponent equal to index equal to $|G|$.\end{remark}

For the remainder of this section let $G=\G$ a noncyclic finite abelian $p$-group with order $p^n$ and $F$ a field with $G$-action so that $F^*$ is $H^1$-trivial.  Let $\Delta'(G)$ be the corresponding algebra with exponent equal to index.  Let $Y=SB(\Delta'(G)^{\otimes p})$, the Severi-Brauer variety of the $p$-th tensor power of $\Delta'(G)$.  Let $\cF(Y)$ be the function field of $Y$, that is, the Severi-Brauer splitting field of the algebra $\Delta'(G)^{\otimes p}$.  By \cite{Saltman-invariant} Theorem 0.5 or \cite{LN} 13.15, $\cF(Y) \cong F(\mpr)^G$.  Here $\mpr$ is a $G$-lattice designed to trivialize the $p$-th power of $[c_2]$.  We define $\mpr$ as follows.  
Let $c_2$ be a 2-cocycle in the class $[c_2]$.  As an abelian group $\mpr \cong \atg \oplus \ig$.  As a $G$-module the $G$-action on $\mpr$ is defined by 
\[g(0,g'-1) = (p\cdot c_2(g,g'),g(g'-1)).\]  
Define $\Delta(G)=\Delta'(G)\otimes_{L_0^G}\cF(Y)$.  Since $F(\atg) \cap F(\mpr)^G=F(\atg)^G$, with the intersection taking place in $F(\mpr)$, the field join of $F(\atg)$ and $F(\mpr)^G$ is $F(\mpr)$ and therefore,
\begin{equation}
\Delta(G)=\Delta'(G)\otimes_{L_0^G}\cF(Y)\cong (F(\mpr)/F(\mpr)^G,z_{\sigma},e(u),e(b)),\label{eqn76}
\end{equation}
where $e(u)$ and $e(b)$ are from $\Delta'(G)$.  The algebra $\Delta(G)$ is an example of the ``generic'' algebras which are treated by Karpenko in \cite{Karpenko2} Definition 3.12.  We recall that definition here.
\begin{definition}[\cite{Karpenko2} Definition 3.12]  Fix a sequence of positive integers $n=n_0>n_1>\cdots > n_m=0$ and a prime $p$.  Let $A$ be a division algebra with $\textrm{ind}(A)=\exp(A)=p^n$.  For each $i=1,2,\ldots,m$ consider the generalized Severi-Brauer variety $Y_i=SB(p^{n_i},A^{\otimes p^i})$, the variety of rank $p^{n_i}$ left ideals in $A^{\otimes p^i}$.  Denote the function field $\cF(Y_1\times \cdots \times Y_m)$ by $\widetilde{F}$ and put $\widetilde{A}=A_{\tilde{F}}$.  Any algebra constructed in this way is a {\bf generic $p$-primary division algebra}.
\end{definition}
Since $\exp(\Delta'(G))=\ind(\Delta'(G))$, the algebra $\Delta(G)=\Delta'(G)_{\mathcal{F}(Y)}$ from \eqnref{eqn76} is an example of a \emph{generic} $p$-primary division algebra with $m=1$, and fixed sequence $n=n_0>n_1=0$ (\cite{Karpenko2} ex. 4.12).  By construction $\Delta(G)$ has $\mathrm{deg}(\Delta(G))=p^n$ and $\exp(\Delta(G))=p$.  We recall Karpenko's calculation of the torsion in $\CH^2$ for certain algebras of exponent $p$ in the next proposition.
\begin{proposition}[\cite{Karpenko2} Proposition 5.1]  Let $A$ be an algebra of prime exponent $p$ and let $X=SB(A)$.  Then the group $\tors \CH^2(X)$ is trivial or (cyclic) of order $p$.  It is trivial if $\textrm{ind}(A)=p$ or $\textrm{ind}(A)\mid 4$.  It is not if $A$ is a ``generic'' division algebra of index $p^n$ and exponent $p$ where $n \geq 2$ in the case of odd $p$ and $n \geq 3$ in the case when $p=2$.
\label{propK}
\end{proposition}
In \cite{Karpenko2} Proposition 5.3, Karpenko proves  if $A$ is a decomposable division algebra of prime exponent, then the group $\CH^2(SB(A))$ is torsion free.  Combining this with \cite{Karpenko2} Proposition 5.1, Karpenko shows that all \emph{generic} algebras of prime exponent $p$ and index $p^n$ are indecomposable, excluding the Albert case ($p=n=2$).  Next we use \cite{Karpenko2} Proposition 5.1 to calculate degeneracy in abelian crossed products.
\begin{corollary}  Let $G=\G$ be a noncyclic abelian $p$-group of order $p^n$, let $F$ be a field with a $G$-action so that $F^*$ is $H^1$-trivial.  Let $\Delta(G)$ be as in \eqnref{eqn76}.  Then $\exp(\Delta(G))=p$ and $\deg(\Delta(G))=p^n$.  If $p\ne 2$ the matrix $e(u)$ is non-degenerate in $F(\mpr)$.  If $p=2$ and $r \geq 3$, the matrix $e(u)$ is not strongly degenerate in $F(\mpr)$.
\label{cor22}
\end{corollary}
\proof Let $X=SB(\Delta(G))$.  $\Delta(G)$ is a ``generic'' division algebra with exponent $p$ and degree $p^n$.  In both cases ($p\ne 2$ or $p=2$ and $r\geq 3$), by \cite{Karpenko2} Proposition 5.1, $\tors \CH^2(X)$ is cyclic of order $p$.  If $e(u)$ is degenerate (resp. strongly degenerate for $p=2$, $n \geq 3$), then $\CH^2(X)$ is torsion free by Proposition \ref{prop7}, a contradiction.\endproof

\begin{remark}  The author has obtained the results of Corollary \ref{cor22} using elementary methods involving studying the lattices defining the field extension $F(\mpr)$ (see \cite{M}).  The proof using lattices is much longer and more involved than the proof presented here, but yields an elementary alternative to using torsion in $\CH^2$ to prove that degeneracy of the matrix will not occur in $\Delta(G)$ for $p\ne 2$.   When $p=2$ the matrix defining $\Delta(G)$ is always degenerate.  This is also proven in \cite{M}, it is essentially a corollary to Albert's result that degree 4, exponent 2 algebras are biquaternion.  The lattice method also yields that for $p=2$ the matrix $e(u)$ of Corollary \ref{cor22} is not strongly degenerate.  In fact, for the case $p=2$, $r=n=3$, a MAGMA program has been written to check for strong degeneracy in $F(\mpr)$, showing that the matrix $e(u)$ is not strongly degenerate in this case.  One can prove that it is enough to check the case $p=2$, $r=n=3$ to cover all cases $p=2$, $r,n \geq 3$.
\end{remark}

\subsection{Examples}
\begin{example}[$p\ne 2$]  Let $p$ be an odd prime, $G$ a noncyclic 
abelian $p$-group and $F$ a field of characteristic $p$.  Let $F^*$ be a $G$-module with trivial action so that $F^*$ is $H^1$-trivial.  Let $\Delta(G)$ be the abelian crossed product as in \eqnref{eqn76} associated to the group $G$ and the field $F$.  That is, $\Delta(G) = \Delta'(G)\otimes \cF(Y)$ where $Y=SB(\Delta'(G)^{\otimes p})$ and $\Delta'(G)$ is the generic crossed product defined by the group $G$.  By Corollary \ref{cor22} the matrix defining the abelian crossed product $\Delta(G)$ is non-degenerate.  Therefore, by Corollary \ref{cor16}, $\cA_{\Delta(G)}$, the generic abelian crossed product associated to $\Delta(G)$, is indecomposable of index $|G|$ and exponent 
the exponent of $G$.  Moreover $\cA_{\Delta(G)}$ remains indecomposable after any prime to $p$ extension.   Note that we may obtain the lowest possible exponent by taking $G$ to be elementary abelian, for then $\cA_{\Delta(G)}$ is indecomposable of index $|G|$ and exponent $p$.  Other examples of indecomposable $p$-algebras can be found in \cite{Jacob} and characteristic independent examples are given in \cite{Karpenko2}.  The algebras in \cite{Jacob} are of index $p^2$ and exponent $p$, $p\ne2$.  \label{example1}
\end{example}
\begin{example}[$p=2$]
Let $F$ be a field of characteristic 2 and let $G = \langle \sigma_1 \rangle \times \langle \sigma_2\rangle\times\langle\sigma_3\rangle$ with $|\sigma_i|=2$, $i=1,2,3$.  Let $\Delta(G)$ be the abelian crossed product given in \eqnref{eqn76} associated to $G$.  By Corollary \ref{cor22} the matrix defining $\Delta(G)$ is not strongly degenerate.  Therefore, by \cite{mck-primetop} Theorem 3.7
, $\cA_{\Delta(G)}$, the generic abelian crossed product associated to $\Delta(G)$, contains no square central elements even after a prime to $2$ extension.  In particular $\cA_{\Delta(G)}$ is indecomposable of exponent 2, index 8 and remains so after any prime to $2$ extension.  Rowen obtained the first example of an indecomposable 2-algebra of exponent 2 and index 8 in \cite{Rowen}. \label{example2}
\end{example}
\begin{remark}In both of these examples note that by \cite{Karpenko2} Proposition 5.3, $\Delta(G)$ itself is indecomposable of exponent $p$ and remains indecomposable after any prime to $p$ extension.  We have proven that the generic abelian crossed products associated to these indecomposable abelian crossed products are also indecomposable.  However, this fact is not used in proving that the generic abelian crossed products are indecomposable.  In fact it is not necessary that the underlying abelian crossed product be indecomposable for the associated generic abelian crossed product to be indecomposable.  The author has construced a decomposable abelian crossed product $p$-algebra with non-degenerat matrix, $p\ne 2$, whose associated generic abelian crossed product is indecomposable.  Also, as mentioned following Corollary \ref{cor22}, the indecomposability of these generic abelian crossed product $p$-algebras is obtainable in a purely algebraic fashion without the use of the Chow group.
\end{remark}
\bibliographystyle{alpha} 
\bibliography{ref.bib}

\end{document}